\documentclass{ws-sd}
\usepackage[T1]{fontenc}
\usepackage[utf8]{inputenc}
\usepackage{newtxtext}
\usepackage{amsmath}
\usepackage{newtxmath,bm}
\usepackage{graphicx}
\usepackage[section]{placeins}
\usepackage{flafter}
\usepackage{algorithm,algpseudocode}
\newcommand{\dd}{\,\mathrm{d}}
\newcommand{\R}{\mathbb{R}}
\newcommand{\EX}{\mathbb{E}}
\newcommand{\PX}{\mathbb{P}}
\newtheorem{assumption}{Assumption}[section]

\begin{document}
% Suppress running headers and the template production header; retain footer page numbers.
\pagestyle{plain}
\renewcommand{\draftnote}{}
% Publisher and journal catchline omitted.
\title{DENSITY EVOLUTION FOR STOCHASTIC DYNAMICAL SYSTEMS\\
DRIVEN BY $\alpha$-STABLE PROCESSES:\\
AN EQUATION-FREE NUMERICAL APPROACH}
\author{YAN LIU}
\address{School of Mathematics and Statistics,\\
Huazhong University of Science and Technology,\\
Wuhan, Hubei, China\\
liuyan1021@hust.edu.cn}
\author{XU SUN\footnote{Corresponding author.}}
\address{School of Mathematics and Statistics,\\
Huazhong University of Science and Technology,\\
Wuhan, Hubei, China\\
xsun@hust.edu.cn}
\maketitle
% Publication dates and AMS classification are not supplied in the source.
% Add publisher-assigned history and author-confirmed classification here when available.
\begin{abstract}
While well-established approaches exist for computing probability densities of stochastic differential equations driven by Brownian motion, computing probability densities for general stochastic differential equations driven by $\alpha$-stable processes remains challenging, and few general-purpose methods are available. A major difficulty is that the corresponding Fokker--Planck equations governing the density are generally unavailable in explicit form. In this paper, we develop an efficient approach to compute probability densities for stochastic differential equations driven by $\alpha$-stable processes without relying on Fokker--Planck equations or any other governing equations. The proposed approach is easy to implement and computationally efficient. Numerical examples demonstrate its effectiveness and indicate an approximately first-order convergence rate with respect to the time step.
\end{abstract}

\section{\label{sec:level1} Introduction}
Probability density contains all statistical information about a stochastic system, and it is one of the central tasks in stochastic dynamics to quantify how the probability density evolves with time. Stochastic differential equations (SDEs) driven by Brownian motions are appropriate models for  stochastic dynamical systems excited by Gaussian white noise. There are well established approaches to compute probability densities associated with SDEs driven by Brownian motions \cite{bally1996}.  These approaches often involve solving  the  corresponding Fokker-Planck equations.

SDEs driven by $\alpha$-stable processes  have found wide applications in a variety of  fields such as  climate dynamics,  finance,  biology,  signal processing \cite{shao2002signal,DiPaola2005}.  Unlike  SDEs driven by Brownian motions,  computing  the probability densities associated with  SDEs driven by non-Gaussian  L\'{e}vy processes remains  challenging.  The main reason lies in the fact that the corresponding Fokker-Planck equations are  often difficult to derive explicitly or solve numerically.

Consider the following SDE driven by  an  $\alpha$-stable process,
\begin{equation}\label{e1.1}
	\left\{
	\begin{aligned}
		& \dd X(t)=f(X(t))\dd t +g(X(t))\dd L(t),\\
		& X(0)=x_0,
	\end{aligned}\right.
\end{equation}
where $X(t)$ is an $\R^d$-valued stochastic process and $L(t)$ is an $\R^d$-valued $\alpha$-stable process ($\ 0<\alpha\leq 2$) defined on the probability space $(\Omega,\mathcal{F},\PX)$, $f:\R^d\to\R^d\times\R^d$ and $g:\R^d\to\R^d$ are the drift and noise coefficients, respectively, and $x_0\in\R^d$ is a given deterministic vector. Note that an  explicit expression for the Fokker-Planck equation of \eqref{e1.1} is generally unavailable, except for the special case  when the coefficient $g$ in  \eqref{e1.1} is a constant \cite{schertzer2001}.  For this case, even though the  Fokker-Planck equation  can be expressed explicitly,  numerical solution remains difficult due to the presence of the fractional Laplace operator.

Let $ C_b(\R^d)$ denote the space of bounded continuous functions on 
$\R^d$. We will need  the following matrix norm:
\begin{equation}
	\|A\|=\sqrt{\sum_{i,j}a_{ij}^2},
\end{equation}
for $A=(a_{ij})$. Throughout the paper, the following assumptions hold:
\begin{assumption}\label{H1}
	The drift coefficient  $f\in C_b(\R^d)$ is  Lipschitz continuous, i.e., there exists a constant $C>0$ such that 
	\begin{equation}
		|f(x)-f(y)|\leq C|x-y|, \forall x,y\in \R^d. 
	\end{equation}
\end{assumption}
\begin{assumption}\label{H2}
	The  coefficient  $g$  satisfies the following condition: there exists a constant $C>0$ such that
	\begin{equation}
		\|g(x)-g(y)\|\leq C|x-y|, \forall x,y\in \R^d
	\end{equation}
	and $g$ satisfies: there exist constants  $ 0<C_1<C_2<\infty$ such that  for all $x\in\R^d, \xi\in\R^d$,
	\begin{equation}
		C_1|\xi|^2\leq \xi^T g(x)^Tg(x)\xi\leq  C_2|\xi|^2.
	\end{equation}
\end{assumption}
Notice that under Assumptions \ref{H1} and \ref{H2},  \eqref{e1.1} has a unique solution, which is a standard result since the coefficients $f$ and $g$ are Lipschitz \cite{applebaum2009}. Under Assumptions \ref{H1} and Assumption \ref{H2}, the solution $X(t)$ to \eqref{e1.1} has a density for any $t>0$, as in \cite{debussche2013,jacod1987}. Note that Assumptions \ref{H1} and \ref{H2} may be  replaced by any sufficient condition. 

In this paper, we introduce   an efficient numerical method for computing the probability density of the solution to  \eqref{e1.1}. The remainder of this paper is organized as follows. In Section \ref{sec:level2}, the key properties of $\alpha$-stable processes are presented. The main equation is derived, and a numerical method to solve it is developed in Section \ref{sec:level3}. Several numerical examples are presented in Section \ref{sec:level4} to illustrate and verify the proposed method.

\section{\label{sec:level2}Notations and  preliminary results}
In this section,  we introduce some notation and preliminary results which will be used in later sections.

% \subsection{Probability Density for a  $\alpha$-stable L\'{e}vy Process}

A nontrivial  one-dimensional $\alpha$-stable L\'{e}vy  process  $L(t)$ is characterized by its characteristic function
\begin{equation}\label{e2.1}
	\EX(e^{izL(t)})=\exp\{t\eta(z)\},
\end{equation}
where 
\begin{equation}\label{e2.2}
	\eta(z)=\left\{
	\begin{aligned}
		&i\gamma z- \sigma^{\alpha}|z|^\alpha\Big(1-i\beta\tan\frac{\pi\alpha}{2}{\rm sgn }\ (z)\Big),\ &&\alpha\neq 1, \\
		&i\gamma z -\sigma |z|\Big(1+i\beta\frac{2}{\pi} {\rm sgn}(z)\log(|z|)\Big), &&\alpha=1.
	\end{aligned}\right.
\end{equation}
The characteristic function   is uniquely determined by these four parameters:
\begin{itemlist}
	\item $\alpha\in(0,2]$ is the characteristic exponent. Note that $\alpha=2$ corresponds to the Gaussian distribution. 
	\item $\beta\in[-1,1]$ is a symmetry parameter.
	\item $\sigma>0$ is a scaling parameter.
	\item $\gamma\in\R$ is a shift parameter.
\end{itemlist}
The $\alpha$-stable process has stationary increments, i.e.,
\begin{equation}\label{e2.3}
	L(t)-L(s)\overset{d}{=}L(t-s)\sim S_\alpha(\sigma(t-s)^{\frac{1}{\alpha}}, \beta,(t-s)\gamma),
\end{equation}
for any $0\leq s<t<\infty$, where $S_\alpha(\sigma, \beta,\gamma)$ stands for an  $\alpha$-stable random variable whose characteristic function is $\exp\{\eta(z)\}$. 

Although the characteristic function of $L(t)$ provides a complete description, there is no closed-form expression for the probability density of $L(t)$. The density can be represented using a series expansion \cite{applebaum2009}. Without loss of generality, assume that at $t=1$, $L(1) \sim S_\alpha(\sigma,\beta,\gamma)$, and denote its probability density function by $p_1^\alpha(x)$. The density for any $t>0$ can be obtained efficiently using the scaling property of stable distributions: 
\begin{equation}\label{e2.4}
	p^\alpha_t(x)=\left\{
	\begin{aligned}
		&t^{-\frac{1}{\alpha}}p_1^\alpha(t^{-\frac{1}{\alpha}}x+(1-t^{\frac{\alpha-1}{\alpha}})\gamma),\ &\alpha\neq 1,\\
		&t^{-1}p_1^\alpha(t^{-1}x-\frac{2}{\pi}\sigma\beta\log(t)),\ &\alpha=1.
	\end{aligned}\right.
\end{equation}
In particular, if $\gamma = 0$ and $\alpha\neq 1$, \eqref{e2.4} simplifies to
\begin{equation}\label{e2.5}
	p^\alpha_t(x) = t^{-1/\alpha} p_1^\alpha(t^{-1/\alpha} x).
\end{equation} 
Throughout this paper,  we assume $p_1^\alpha(x)$ is known.  In fact,  for an $\alpha$-stable random variable $S_\alpha(\sigma,\beta,\gamma)$, such as $L(1)$,  there are well established approaches to compute its probability density \cite{janicki1993}.

\section{\label{sec:level3}Numerical Method}
In this section, we present the numerical framework for computing the probability density of the solution to \eqref{e1.1}. Note that in this work, we only give the proof for $d=1$, while numerical results are also presented for multi-dimensional systems.
\subsection{ Euler-Maruyama Scheme}\label{section3.1}
Let \(\Delta t=T/N\) and define the uniform time grid \(t_k=k\Delta t\), \(k=0,1,\ldots,N\). Applying the Euler scheme to \eqref{e1.1} gives 
\begin{equation}\label{e2.6}
	\left\{\begin{aligned}
		&\hat{X}(t_{k+1})=\hat{X}(t_{k})+f(\hat{X}(t_{k}))\Delta t+g(\hat{X}(t_{k})) \Delta L(t_k),\\
		&\hat{X}(t_0)=x_0,
	\end{aligned}
	\right.
\end{equation}
for $k=0,1,\cdots,N-1$, where  $ \Delta L(t_k)= L(t_{k+1})-  L(t_k)$.  By the stationary increment property in \eqref{e2.3},   $ \Delta L(t_k)$ has the same distribution as $L(\Delta t)$, i.e., the probability density of $\Delta L(t_k)$  can be expressed  by using $p^\alpha_{\Delta t} (x)$ defined in \eqref{e2.4}.

Throughout the paper,  we use $p$ and $\hat p$ to denote the transition density associated with \eqref{e1.1} and \eqref{e2.6}, respectively.   Let  $\hat{p}(x,t_{k+1}|y,t_k)$ represent the  transition density of $\hat X(t_{k+1})$ at $x$ under the condition $\hat X(t_k)=y$,  while  ${p} (x,t|y,s)$  denote the density of $ X(t)$ at $x$ under the condition $X(s)=y$. It follows from \eqref{e2.6} that 
\begin{equation}\label{e2.7}
	\begin{aligned}
			&\hat{p}(x_{k}, t_{k} | x_{k-1}, t_{k-1})\\
			 =& \frac{1}{|g(x_{k-1})|}
		p_{\Delta t}^\alpha\left(\frac{x_k-x_{k-1}-f(x_{k-1})\Delta t}{g(x_{k-1})}\right).
	\end{aligned}
\end{equation}
Substituting \eqref{e2.4} into \eqref{e2.7}, we obtain that  for $\alpha\neq 1$,
\begin{equation}\label{e2.8}
	\begin{aligned}
		&\hat{p}(x_{k}, t_{k} \mid x_{k-1}, t_{k-1})\\
	    =&
		\frac{1}{|g(x_{k-1})|(\Delta t)^{1/\alpha}}
		p_1^\alpha
		\Bigg(
		\frac{x_k-x_{k-1}-f(x_{k-1})\Delta t}
		{g(x_{k-1})(\Delta t)^{1/\alpha}}\\
		&
		+\left[1-(\Delta t)^{(\alpha-1)/\alpha}\right]\gamma
		\Bigg),
	\end{aligned}
\end{equation}
and for $\alpha = 1$,
\begin{equation}\label{e2.9}
	\begin{aligned}
		&\hat{p}(x_{k}, t_{k} \mid x_{k-1}, t_{k-1})\\
		=& \frac{1}{|g(x_{k-1})|\Delta t}
		{p_1^\alpha}\Bigg(
		\frac{x_k-x_{k-1}-f(x_{k-1})\Delta t}
		{g(x_{k-1})\Delta t}
		-\frac{2}{\pi}\sigma\beta\log(\Delta t)
		\Bigg).
	\end{aligned}
\end{equation}

Denote $\hat{p}_{\mathcal{A}}(x,t_{k})$ as the probability density of $\hat X(t_k)$ defined in \eqref{e2.6},  and ${p}_{\mathcal{A}}(x, t_{k})$ the probability density of $X(t)$ defined in \eqref{e1.1}.  Here,  the subscript '$\mathcal{A}$' is used to explicitly indicate the initial condition.  $\hat{p}_{\mathcal{A}}(x,t_{k})$ and ${p}_{\mathcal{A}}(x, t_{k})$ can be expressed in terms of the conditional density as
\begin{align}
	\hat{p}_{\mathcal{A}}(x,t_{k})&=\hat{p} (x,  t_k|x_0,0),\label{e2.10}\\
	{p}_{\mathcal{A}}(x,t)&={p} (x,  t|x_0,0),\label{e2.11}
\end{align}
respectively.

In this work, we assume that ${\hat p}_{\mathcal{A}}(x,t_k)$ converges to  ${p}_{\mathcal{A}}(x,t_k)$ as $\Delta t \rightarrow 0 $ for all $x\in\R$. Sufficient conditions for density convergence have been established \cite{knopova2018,konakov2011}. In particular, when $L(t)$ is a symmetric $\alpha$-stable process with $\alpha\in(1,2)$, under suitable regularity and non-degeneracy assumptions on the coefficients, the transition probability density ${\hat p}_{\mathcal{A}}$ associated with the Euler approximation satisfies the following error estimate \cite{konakov2011}:
\begin{equation}\label{e2.12}
	\left|{\hat p}_{\mathcal{A}}(x,t_k)-{p}_{\mathcal{A}}(x,t_k)\right|
	\leq C\Delta t.
\end{equation}
In the general setting considered in this work, particularly for non-symmetric $\alpha$-stable processes,  related convergence results for Euler
approximations at the level of the underlying stochastic processes  available in the literature \cite{janicki1996}.  However, establishing the convergence of the corresponding probability densities for general $\alpha$-stable noise, including both symmetric and non-symmetric cases over the full range $\alpha\in(0,2]$, is beyond the scope of this paper and is left for future work. The convergence behavior with respect to $\Delta t$ is examined numerically in Example \ref{example4.1}.

\subsection{Numerical Method}\label{section3.3} 
In this section, we develop an efficient numerical algorithm to compute  ${\hat p}_{\mathcal{A}}(x,t_{k})$ .

It follows from \eqref{e2.10} that for $k=2, \dots,N$,
\begin{equation}\label{e3.1}
	\begin{aligned}
		&\hat{p}_{\mathcal{A}}(x,t_k) =\hat{p}(x,t_k|x_0,0)\\
		=& \int_{\R}  \hat{p}(x,t_k;x_{k-1},t_{k-1}|x_0,0)\;\dd x_{k-1}\\
		=& \int_{\R} \hat{p}(x,t_k|x_{k-1},t_{k-1}; x_0,0)\hat{p}(x_{k-1},t_{k-1}|x_0,0)\; \dd x_{k-1}.
	\end{aligned}
\end{equation}
It is straightforward to check that $\hat X(t_k)$  defined in \eqref{e2.6} has  Markovian and time-homogeneous properties. Thus, for $j=1,2, \cdots,N$,
\begin{equation}\label{e3.2}
	\begin{aligned}
			&\hat p(x_j,t_j|x_{j-1},t_{j-1};\cdots; x_0,0)\\
			=&\hat p(x_j,t_j|x_{j-1},t_{j-1})\\
			=&\hat p(x_j,\Delta t|x_{j-1},0).
	\end{aligned}
\end{equation}
Substituting \eqref{e3.2} into  \eqref{e3.1},  we obtain
\begin{equation}\label{e3.3}
	\begin{aligned}
		\hat p_{\mathcal{A}}(x,t_k)=&\int_{\R} \hat{p}(x,t_k|x_{k-1},t_{k-1})\hat{p}(x_{k-1},t_{k-1}|x_0,0)\; \dd x_{k-1}\\
		=&\int_{\R} \hat{p}(x,\Delta t |x_{k-1},0)\hat p_{\mathcal{A}}(x_{k-1},t_{k-1})\; \dd x_{k-1}.
	\end{aligned}
\end{equation}
For $k=1$, the transition probability density of $\hat{X}(t_1)$ is  
\begin{align}\label{e3.4}
	\hat {p}_{\mathcal{A}}(x,t_1)=\hat{p}(x,t_1|x_0,0),
\end{align}
which has an explicit expression  in \eqref{e2.7}.

%It follows from Eq.~\eqref{e3.3} that
%\begin{equation}\label{e3.5}
%	\begin{aligned}
%		&\hat p_{\mathcal{A}}(x,t_k)\\
%		= &\int_{\R} \hat{p}(x,\Delta t|x_{k-1},0)\hat p_{\mathcal{A}}(x_{k-1},t_{k-1})\;\dd x_{k-1}.
%	\end{aligned}
%\end{equation}
Eq.~\eqref{e3.3} can be approximated as
\begin{equation}\label{e3.6}
	\begin{aligned}
		%&\hat p_{\mathcal{A}}(x,t_k)\\
		&\int_{-D}^D \hat{p}(x,\Delta t|x_{k-1},0)\hat p_{\mathcal{A}}(x_{k-1},t_{k-1})\;\dd x_{k-1},
	\end{aligned}
\end{equation} 
where $D$ is a sufficiently large constant.

Discretize the interval $[-D,D]$ with $M+1$ equidistant grid points $z_0,z_1,\cdots, z_M$ such that $z_k=-D+k\Delta x$ for $k=0,1,2,\cdots,M$ and $\Delta x=\frac{2D}{M}$. Discretizing \eqref{e3.6} by the trapezoidal rule yields 
\begin{equation}\label{e3.7}
	\begin{aligned}
		&\hat p_{\mathcal{A}}(x,t_k)
		%=&\sum_{i_d=0}^{M-1}\overbrace{\int_{-D}^D\cdots\int_{-D}^D}^{m-1}\int_{z_{i_d}}^{z_{i_d+1}} \hat{p}(x^1,x^2,\cdots x_k^d,t_k|x_{k-1}^1,x_{k-1}^2,\cdots,x_{k-1}^d,t_{k-1})\hat p_{\mathcal{A}}(x_{k-1}^1,x_{k-1}^2,\cdots,x_{k-1}^d,t_{k-1})\;\prod_{j=1}^d\dd x_{k-1}^j\\
		%=&\sum_{i_d=0}^{M}\overbrace{\int_{-D}^D\cdots\int_{-D}^D}^{m-1} \omega_{i_d}\hat{p}(x^1,x^2,\cdots x_k^d,t_k|x_{k-1}^1,x_{k-1}^2,\cdots z_{i_d},t_{k-1})\hat p_{\mathcal{A}}(x_{k-1}^1,x_{k-1}^2,\cdots,z_{i_d},t_{k-1})\Delta t+\mathcal{O}(\D t^3)\;\prod_{j=1}^{d-1}\dd x_{k-1}^j\\
		=\Delta x\sum_{i=0}^{M} \omega_{i}\hat{p}(x,\Delta t|z_{i},0)\hat p_{\mathcal{A}}(z_{i},t_{k-1})+\mathcal{O}(\Delta x^2),
	\end{aligned}
\end{equation}
where $\hat{p}(x,\Delta t|z_{i},0)$ can be obtained from Eqs.~\eqref{e2.8} and \eqref{e2.9}, and  the constants $\omega_{i}$ satisfy
\begin{equation}\label{e3.8}
	\omega_{i}=\left\{\begin{aligned}
		&\frac{1}{2},\quad \text{for } i=0 \text{ or } M,\\
		&1,\quad  \text{ for } i=1,2,\cdots,M-1.
	\end{aligned}\right.
\end{equation}

Omitting the small error term in \eqref{e3.7}, we obtain $\tilde{p}$, defined as
\begin{equation}\label{e3.9}
	\tilde{p}_{\mathcal{A}}(x,t_1)=\hat {p}_{\mathcal{A}}(x,t_1)=\hat{p}(x,t_1|x_0,0),
\end{equation}
and for $k=2,3,\cdots,N$,
\begin{equation}\label{e3.10}
	\begin{aligned}
		&\tilde{p}_{\mathcal{A}}(x,t_k)
		=\Delta x\sum_{i=0}^{M} \omega_{i}\hat{p}(x,\Delta t|z_{i},0)\tilde{p}_{\mathcal{A}}(z_{i},t_{k-1}).
	\end{aligned}
\end{equation}
Comparing \eqref{e3.10}  with \eqref{e3.7}, we obtain   
\begin{equation}\label{e3.11}
	\hat{p}_{\mathcal{A}}(x,t_k)=\tilde{p}_{\mathcal{A}}(x,t_k)+\mathcal{O}(\Delta x^2).
\end{equation}

In summary, the transition probability ${p}_{\mathcal{A}}(x,t_k)$ of $X(t)$ defined in \eqref{e1.1} can  be approximated  by $\tilde{p}_{\mathcal{A}}(x,t_k)$  in \eqref{e3.9} and \eqref{e3.10} as $\Delta t\rightarrow 0$ and $\Delta x \rightarrow 0$. For   $d=1$, 
\begin{equation}\label{e3.12}
	\begin{aligned}
		\tilde {p}_{\mathcal{A}}(z_i, t_1)= \hat{p}(z_i,t_1|x_0,t_0) \quad\text{for $k=1$ and  $i=0, 1, \dots,M$.}
	\end{aligned}
\end{equation}
and
\begin{align}\label{e3.13}
	\tilde {p}_{\mathcal{A}}(z_i, t_k)= \sum_{j=0}^M \omega_j\hat{p} (z_i,\Delta t|z_j,0)\tilde {p}_{\mathcal{A}}(z_j, t_{k-1})\Delta x,
\end{align}
for $k=2,3,\dots,N$ and $i=0, 1, \dots,M$.
Denote $\mathbf{p}_k$ and $A$ as
\begin{align}\label{e3.14}
	\mathbf{p}_k=\left ( \begin{matrix}
		\tilde {p}_{\mathcal{A}}(z_0, t_k)\\
		\tilde {p}_{\mathcal{A}}(z_1, t_k) \\
		\vdots\\
		\tilde {p}_{\mathcal{A}}(z_M, t_k)
	\end{matrix} \right )
\end{align}
and 
\begin{align}\label{e3.15}
	A=(a_{ij})_{(M+1)\times (M+1)} \quad\text{with}\quad a_{ij}=\omega_j\hat{p} (z_i,\Delta t|z_j,0)\Delta x, 
\end{align}
respectively. Then, Eq.~\eqref{e3.13}  can be written as  
\begin{align}\label{e3.16}
	&\mathbf{ p}_{k}=A\,{\mathbf{p}}_{k-1} \quad\text{for $k=2,3,\dots,N$,} 
\end{align}
with initial value $\mathbf{p}_1$ as
\begin{align}\label{e3.17}
	\mathbf{p}_1=\left ( \begin{matrix}
		\tilde {p}(z_0,\Delta t|x_0,0)\\
		\tilde {p}(z_1,\Delta t|x_0,0) \\
		\vdots\\
		\tilde {p}(z_M,\Delta t|x_0,0)
	\end{matrix} \right ).
\end{align}
Note that both the matrix $A$ and the initial value $\mathbf{p}_1$ are constant,  which can be obtained by substituting \eqref{e2.4} into \eqref{e3.15} and \eqref{e3.17}, respectively.  Once $A$ and $\mathbf{p}_1$ are given,  then $\mathbf{p}_k$,   which is the numerical solution to the probability density for \eqref{e1.1},  can be computed recursively by using \eqref{e3.16}.

For systems with $d>1$, we obtain the following results, which are similar to \eqref{e3.12} and \eqref{e3.13}, for $i_1,i_2,\cdots,i_d=0,1,\cdots,M$,
\begin{equation}
	\begin{aligned}
		&\tilde{p}_{\mathcal{A}}(z_{i_1},z_{i_2},\cdots,z_{i_d},t_1)\\
		=&\hat {p}_{\mathcal{A}}(z_{i_1},z_{i_2},\cdots,z_{i_d},t_1)\\
		=&\hat{p}(z_{i_1},z_{i_2},\cdots,z_{i_d},t_1|x_0,0),
	\end{aligned}
\end{equation} 
and for $k=2,3,\cdots,N$,
\begin{equation}
	\begin{aligned}
		&\tilde{p}_{\mathcal{A}}(z_{i_1},z_{i_2},\cdots z_{i_d},t_k)\\
		=&(\Delta x)^d\sum_{j_1,\cdots,j_d=0}^{M} \Big(\prod_{l=1}^d\omega_{j_l}\Big) \hat{p}(z_{i_1},z_{i_2},\cdots z_{i_d},\Delta t|z_{j_1},z_{j_2},\cdots z_{j_d},0)\\
		&\times\tilde{p}_{\mathcal{A}}(z_{j_1},z_{j_2},\cdots z_{j_d},t_{k-1}),           
		\end{aligned}
\end{equation}
where the constants $\omega_{j_l}$ satisfy
\begin{equation}
	\omega_{j_l}=\left\{\begin{aligned}
		&\frac{1}{2},\quad \text{for } j_l=0 \text{ or } M,\\
		&1,\quad  \text{ for } j_l=1,2,\cdots,M-1.
	\end{aligned}\right.
\end{equation}

The steps to implement the proposed algorithm,  which consists of Eqs.~\eqref{e2.4}, \eqref{e3.15},  \eqref{e3.16} and \eqref{e3.17},  are summarized in Algorithm \ref{algorithm}.

\begin{algorithm}[htbp]
	\caption{Compute the probability density of the solution for SDE \eqref{e1.1} with initial value $x_0$ } 
	\label{algorithm}
	
	\begin{algorithmic}[1]
		
		\Require $T, D, \Delta t,\Delta x, f ,g ,\alpha,\beta,\sigma,\gamma$. Here, $f$ and $g$ are functions used in \eqref{e1.1}, $\alpha,\beta,\sigma,\gamma$ are parameters for $L(t)$

		\Ensure $\tilde{p}(x,T)$, which is defined in \eqref{e3.9} and \eqref{e3.10}. $\tilde{p}(x,T)$ stores all the numerical solutions of the probability density values at all the discrete nodes
		
		\State Compute the number of time nodes $ N=\frac{T}{\Delta t}$ and the number of space nodes $M= \frac{2D}{\Delta x}$

		\For {$k=0\  \text{to}\ N$}
		\State Set $t_k=k\Delta t$
		\EndFor
		
		\For {$k=0\  \text{to}\ M$}
		\State Set $z_k=-D+k\Delta x$
		\EndFor
		
		\For {$i=0\ \text{to}\ M$}
		\For{$j=0\ \text{to}\ M$}
		\State Compute $ a_{(i+1)(j+1)}=\omega_j\hat{p} (z_i,\Delta t|z_j,0)\Delta x$ by \eqref{e2.4}
		\EndFor
		\EndFor
		
		\For {$i=0\ \text{to}\ M$}
		\State Compute $	(\mathbf{p}_1)_i=\hat{p}(z_{i},\Delta t|x_0,0)$ by \eqref{e2.4} 
		\EndFor
		
		\State    Get $A$  with $a_{ij}$ obtained
		
		\State Get $\mathbf{p}_N$ by \eqref{e3.16}
		
		\State Output $\mathbf{p}_N$

	\end{algorithmic}
\end{algorithm}

%%%%%%%%%%%%%%%%%%%%%%%%%%%%%%%%%%%%%%%%%%%%%%%%%%%%%%%%

\section{\label{sec:level4}Examples}
In this section, several examples are presented to illustrate and verify the numerical method. 
\begin{example}\label{example4.1}\rm
	Consider the following SDE:
	\begin{equation}\label{e4.1}
		\left\{
		\begin{aligned}
			& \dd X(t)=-X(t)\dd t+\dd L(t),\\
			& X(0)=x_0=0,
		\end{aligned}\right.
	\end{equation}
	where $L(t)$ is a $\R$-valued symmetric $\alpha$-stable process with $\sigma=1, \beta=0, \gamma=0$.
	Eq.~\eqref{e4.1} is a special case of \eqref{e1.1} with $f(x)=-x$ and $g(x)=1$. The solution to \eqref{e4.1} is 
	\begin{equation}
		X(t)=\int_0^t e^{s-t}\dd L(s).
	\end{equation}
	The characteristic function of $X(t)$ is  given by \cite{sato1990}
	\begin{equation}
		\begin{aligned}
				\EX[e^{izX(t)}]=&\exp\Big\{-|z|^\alpha\int_0^t e^{-(t-s)\alpha}\dd s\Big\}\\
				=&\exp\Big\{-|z|^\alpha\frac{1-e^{-\alpha t}}{\alpha}\Big\}.
		\end{aligned}
	\end{equation}
 Therefore, for any $t>0$, the distribution of $X(t)$ is
	\begin{equation}
		X(t)\sim S_\alpha((\frac{1-e^{-\alpha t}}{\alpha})^\frac{1}{\alpha},0,0).
	\end{equation}
	Consequently, the probability density ${p}_{\mathcal{A}}(x,t)$ can be obtained analytically.

	We compute the probability density of $X(T)$ at $T=1$ for $\alpha=0.8$ and $\alpha=1.5$. In the numerical computation, the time step and the space step are taken as $\Delta t = 1/20$ and $\Delta x = 1/100$, respectively, with the spatial domain $[-30,30]$. The  probability densities  of $X(t)$ obtained by the numerical method in Section \ref{section3.3} are compared with the analytical solutions and  Monte Carlo simulations based on $10^5$ sample paths.
	
\begin{figure}[!htbp]
	\centering
	\includegraphics[width=0.80\linewidth]{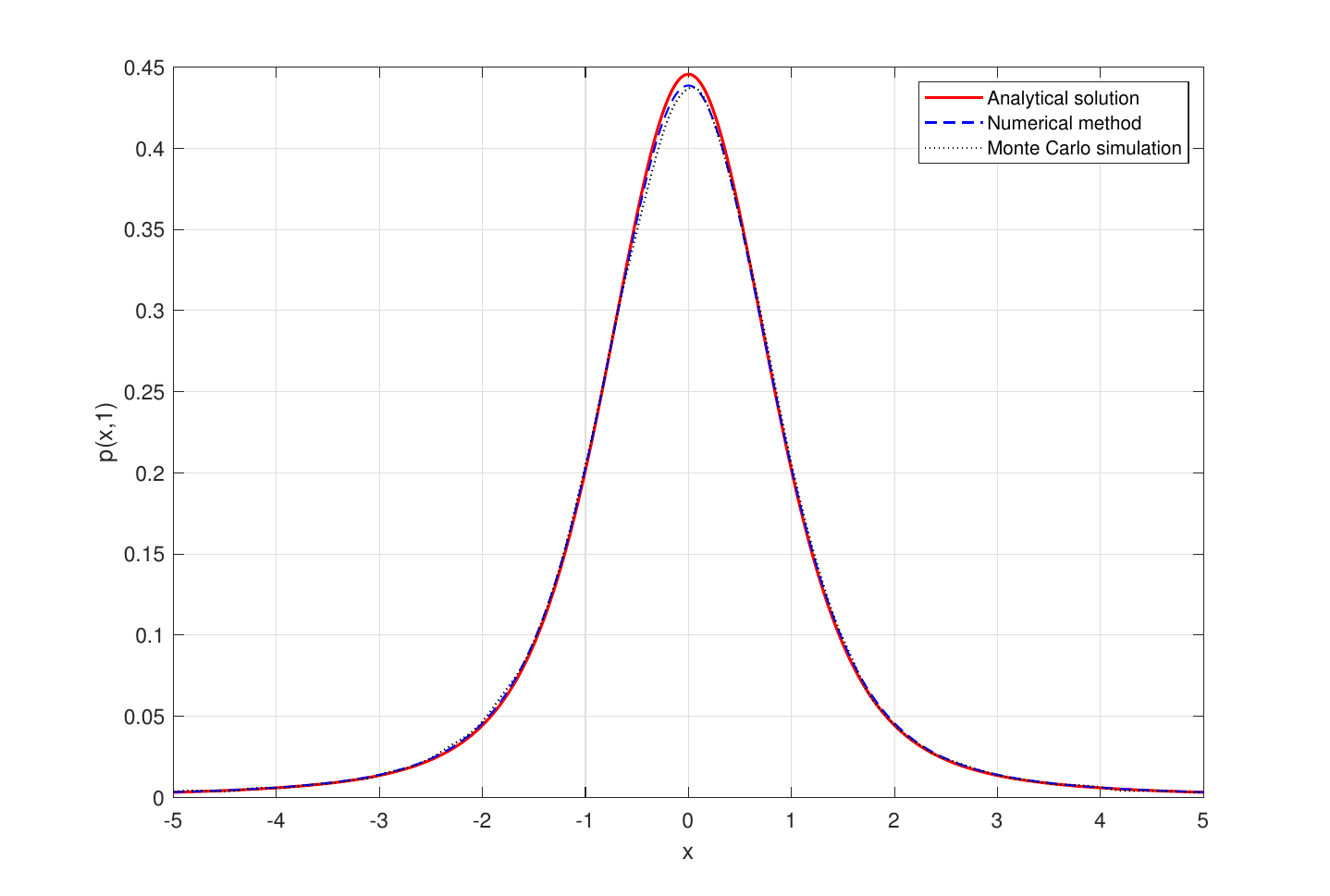}% Here is how to import EPS art
\centering
		\caption{Probability density of the solution to \eqref{e4.1} at $T=1$ for $\alpha=1.5$ ($\Delta t=1/20$, $\Delta x=1/100$): comparison of the numerical method with the analytical solution and Monte Carlo simulation.}
		\label{fig:1}
\end{figure}
	
	\begin{figure}[!htbp]
		\centering
		\includegraphics[width=0.80\linewidth]{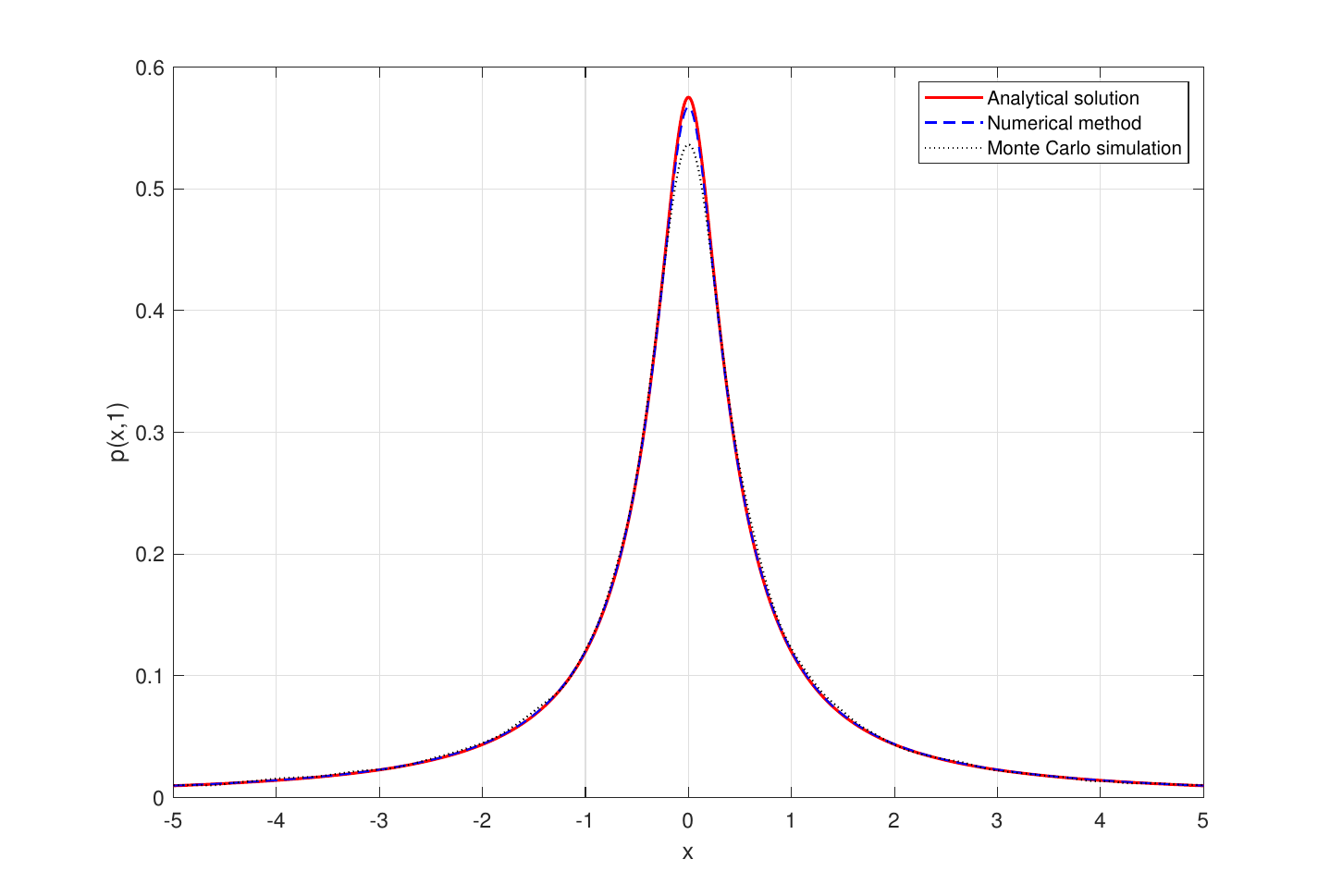}
		\caption{Probability density of the solution to \eqref{e4.1} at $T=1$ for $\alpha=0.8$ ($\Delta t=1/20$, $\Delta x=1/100$): comparison of the numerical method with the analytical solution and Monte Carlo simulation.}
		\label{fig:2}
	\end{figure}

	For $\alpha=1.5$, Fig.~\ref{fig:1} shows that the probability density obtained by the numerical method is nearly indistinguishable from the analytical solution and agrees well with the Monte Carlo result. For the case $\alpha=0.8<1$, Fig.~\ref{fig:2} shows that the numerical result again agrees well with both the analytical solution and Monte Carlo estimate. Compared with the case $\alpha=1.5$, the density for $\alpha=0.8$ exhibits more pronounced heavy-tailed behavior, which is well captured by the numerical method. These results demonstrate that the proposed method accurately computes the probability density for stable noise with both $\alpha>1$ and $\alpha<1$.

	To assess the convergence of the numerical method, we compute the mean-square error (MSE) and the root-mean-square error (RMSE) between the numerical and the analytical densities for the time steps $\Delta t = 1/5,\,1/10,\,1/20$, with $\Delta x=1/100$. The experimental convergence order is estimated by
	\begin{equation}
		\text{Order} = \log_2\Big(\frac{\mathrm{RMSE}(\Delta t)}{\mathrm{RMSE}(\Delta t/2)}\Big).
	\end{equation}
	Tables \ref{tab:conv15} and \ref{tab:conv08} list the errors and the corresponding convergence orders for $\alpha=1.5$ and $\alpha=0.8$, respectively. In both cases the errors decrease essentially linearly with $\Delta t$, and the estimated orders are close to one ( $1.02$ and $1.01$ for $\alpha=1.5$; $1.02$ and $0.90$ for $\alpha=0.8$). These results indicate an approximately first-order convergence rate with respect to the time step $\Delta t$. For $\alpha=1.5$, this numerical convergence rate is consistent with the error estimate in \eqref{e2.12}.

	\begin{table}[htbp]
		\centering
		\tbl{Convergence of the numerical density for $\alpha=1.5$ ($D=30$, $\Delta x=1/100$).\label{tab:conv15}}
{\begin{tabular}{@{}ccccc@{}}
\toprule
			$\Delta t$ & $N$ & MSE & RMSE & Order \\
            \colrule
			$1/5$  & 5  & $1.206\times10^{-5}$ & $3.473\times10^{-3}$ & --- \\
			$1/10$ & 10 & $2.913\times10^{-6}$ & $1.707\times10^{-3}$ & 1.025 \\
			$1/20$ & 20 & $7.183\times10^{-7}$ & $8.475\times10^{-4}$ & 1.010 \\
		\botrule
\end{tabular}}
	
	\end{table}
	\begin{table}[htbp]
		\centering
		\tbl{Convergence of the numerical density for $\alpha=0.8$ ($D=30$, $\Delta x=1/100$).\label{tab:conv08}}
{\begin{tabular}{@{}ccccc@{}}
\toprule
			$\Delta t$ & $N$ & MSE & RMSE & Order \\
			\colrule
			$1/5$  & 5  & $6.418\times10^{-6}$ & $2.533\times10^{-3}$ & --- \\
			$1/10$ & 10 & $1.557\times10^{-6}$ & $1.248\times10^{-3}$ & 1.022 \\
			$1/20$ & 20 & $4.499\times10^{-7}$ & $6.708\times10^{-4}$ & 0.896 \\
		\botrule
\end{tabular}}
	
	\end{table}
% Keep all figures and tables within this example.
\FloatBarrier
\end{example}

%%%%%%%%%%%%%%%%%%%%%%%%%%%%%%%%%%%%%%%%%%%%%%%%%%%%%%

 \begin{example}
	\rm
	Consider the following stochastic dynamical system in $\mathbb{R}^2$ driven by independent symmetric $\alpha$-stable processes:
	\begin{equation}
		\begin{cases}
			\dd X_1(t) = \big(-X_1(t)+\sin X_2(t)\big)\,\dd t + \dd L_1(t),\\[2mm]
		    \dd X_2(t) = \big(-X_2(t)+\cos X_1(t)\big)\,\dd t + \dd L_2(t),
		\end{cases}
		\label{ee42}
	\end{equation}
	where $L_1$ and $L_2$ are independent $\mathbb{R}$-valued symmetric $\alpha$-stable processes with $\sigma=1$, $\beta=0$, $\gamma=0$ and $X(0) = x_0 = (0,0)^T$. Unlike the linear system in Example~\ref{example4.1}, a closed-form expression for the probability density is not readily available due to the nonlinear coupling terms  $\sin X_2$ and $\cos X_1$. The proposed method is therefore validated by comparing the computed density with Monte Carlo simulations.

	\begin{figure}[!htbp]
		\centering
		\includegraphics[width=0.72\linewidth]{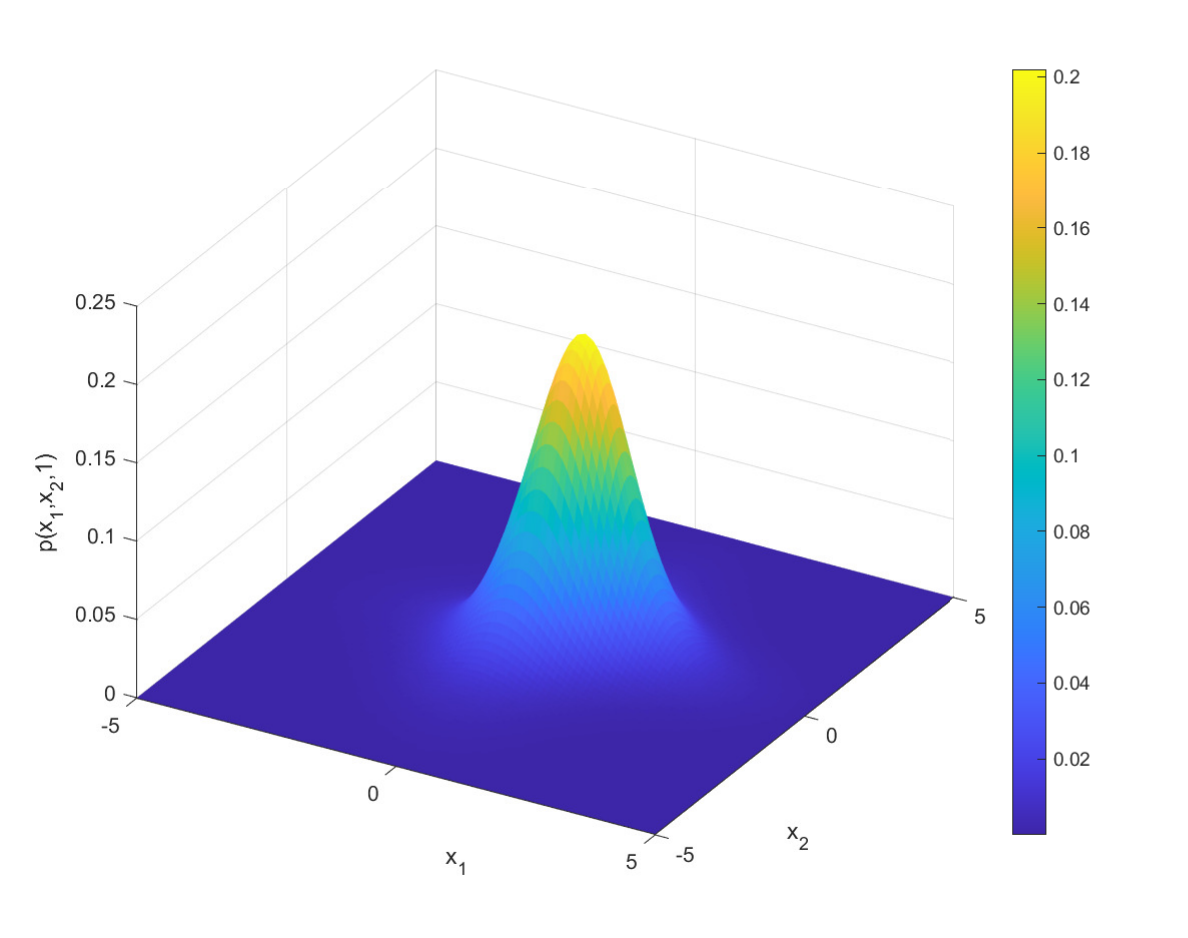}
		\caption{Joint probability density $p(x_1,x_2,1)$ of the solution to \eqref{ee42} computed by the proposed method for $\alpha=1.5$, with $\Delta t=1/20$ and $\Delta x=1/10$.}
		\label{fig:3}
	\end{figure}
	
	In the numerical computation, we take $\alpha=1.5$, $T=1$,  the time step $\Delta t = 1/20$, the space step $\Delta x = 1/10$, and the spatial domain $[-5,5]^2$.  Fig.~\ref{fig:3} shows the joint probability density $p(x_1,x_2,1)$ computed by the numerical method. To validate the numerical results, Monte Carlo simulations are performed using $5\times10^5$ independent sample paths with the same time step $\Delta t=1/20$. Figs.~\ref{fig:4} and \ref{fig:5} compare the marginal probability densities of $X_1$ and $X_2$, respectively, obtained by the numerical method and Monte Carlo simulation. The marginal densities obtained by the two methods are in good agreement, with only slight differences near their peaks. These results demonstrate the applicability and accuracy of the proposed approach for two-dimensional nonlinear stochastic systems driven by $\alpha$-stable processes.

	\begin{figure}[!htbp]
	\centering
	\includegraphics[width=0.80\linewidth]{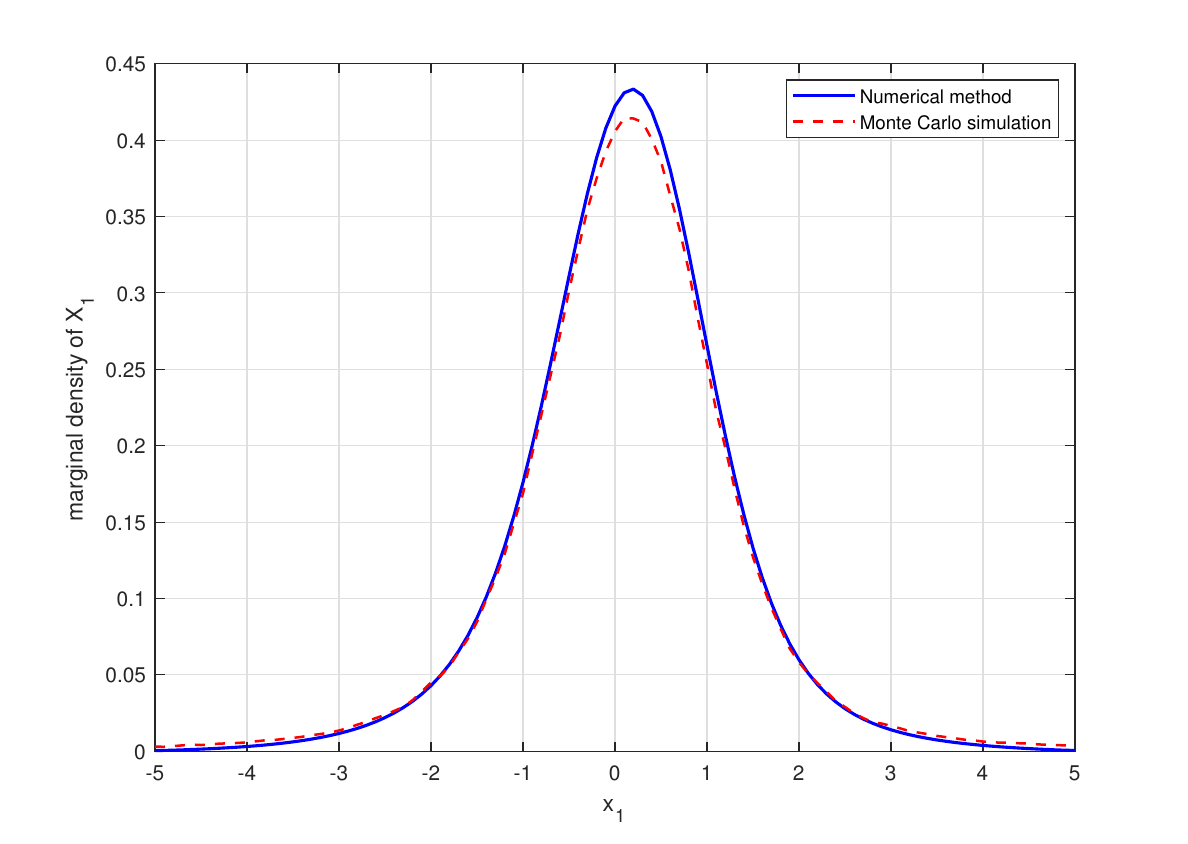}
	\caption{Marginal probability density of $X_1$ at $T=1$ for $\alpha=1.5$: comparison of the proposed method with the Monte Carlo simulation.}
	\label{fig:4}
\end{figure}

	\begin{figure}[!htbp]
	\centering
	\includegraphics[width=0.80\linewidth]{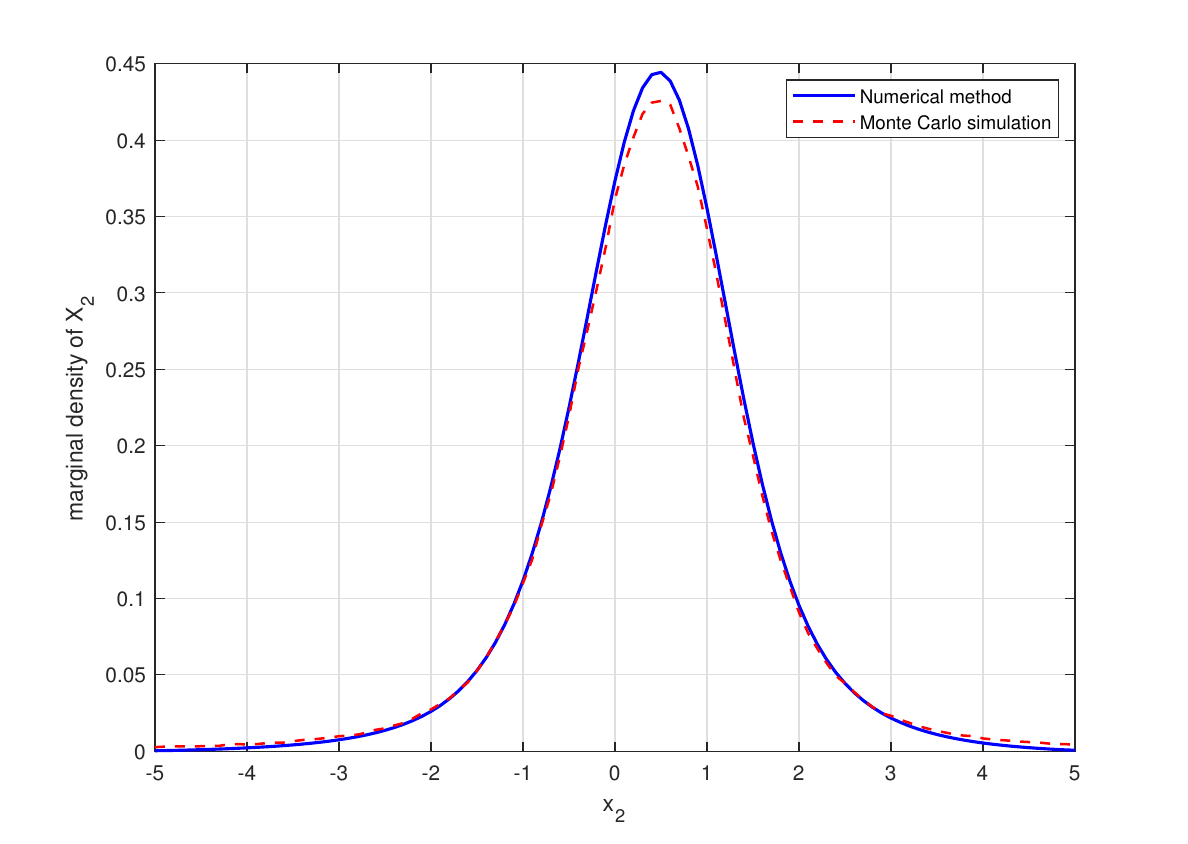}
	\caption{Marginal probability density of $X_2$ at $T=1$ for $\alpha=1.5$: comparison of the proposed method with the Monte Carlo simulation.}
	\label{fig:5}
\end{figure}
%	

% Keep all figures and tables within this example.
\FloatBarrier
\end{example}

\section{\label{sec:level5}Conclusion}
In this paper, we develop an efficient and easy-to-implement numerical method for computing the probability densities of solutions to stochastic differential equations driven by $\alpha$-stable processes, without relying on the derivation or numerical solution of the corresponding Fokker--Planck equations. The proposed method is applicable to both symmetric and non-symmetric $\alpha$-stable processes over the entire range $\alpha\in(0,2]$. Numerical examples demonstrate the accuracy of the proposed method through comparisons with analytical solutions and Monte Carlo simulations, and indicate an approximately first-order convergence rate with respect to the time step. The applicability of the method to multi-dimensional nonlinear systems is also illustrated numerically. Future work will focus on extending the proposed framework to more general stochastic systems, including degenerate systems, with particular interest in applications to stochastic oscillators driven by $\alpha$-stable noise.

%%%%%%%%%%%%%%%%%%%%%%%%%%%%%%%

%\begin{acknowledgments}
%We wish to acknowledge the support of the author community in using
%REV\TeX{}, offering suggestions and encouragement, testing new versions,
%\dots.
%\end{acknowledgments}

%\section*{Declaration}
%The authors  have no conflicts to disclose.

%\section*{Data Availability Statement}

%Data sharing is not applicable to this article as no new data were created or analyzed in this study.

% References are generated from ref.bib with BibTeX.
\FloatBarrier
\bibliographystyle{sd-numeric}
\bibliography{ref}

\end{document}